\documentclass{amsart}

\usepackage[latin1]{inputenc}
\usepackage{comment}
\usepackage{qtree,array,youngtab}

\usepackage{color, graphicx}
\usepackage{soul}
\usepackage[T1]{fontenc}

\definecolor{darkgreen}{rgb}{0.,0.5,0.}

 \newtheorem{thm}{Theorem}[section]
 
 \newtheorem{lem}[thm]{Lemma}
 \newtheorem{ques}[thm]{Question}
\title{Partitions with the same hook multiset}
\author{Huan Xiong}
\address{I-Math, Universit$\ddot{a}$t
Z$\ddot{u}$rich, Winterthurerstrasse 190, Z$\ddot{u}$rich 8057,
Switzerland} \email{huan.xiong@math.uzh.ch}



\begin{document}
\begin{abstract}
It is well-known that two conjugate partitions have the same hook
multiset. But two different partitions with the same hook multiset
may not be conjugate to each other. In $1977$, Herman and Chung
proposed the following question: What are the necessary and
sufficient conditions for partitions to be determined by their hook
multisets up to conjugation? In this paper, we will answer this
question by giving a criterion to determine whether two different
partitions with the same hook multiset are conjugate to each other.
\end{abstract}

%%% ----------------------------------------------------------------------
\maketitle
%%% ----------------------------------------------------------------------
%\tableofcontents

\section{Introduction}
Hook lengths of partitions are very useful in the study of number
theory, combinatorics and representation theory.  A \emph{partition}
is a finite sequence of positive integers $\lambda = (\lambda_1,
\lambda_2, \ldots, \lambda_m)$ such that $\lambda_1\geq
\lambda_2\geq \cdots \geq \lambda_m$. A partition $\lambda$ could be
visualized by its \emph{Young diagram}, which is a finite collection
of boxes arranged in left-justified such that there are exactly
$\lambda_i$ boxes in the $i$-th row. For every box in the Young
diagram of $\lambda,$ we can associate its \emph{hook length}, which
is the number of boxes in the same row to the right of it, in the
same column below it, or the box itself. We denote the hook length
of $(i, j)$-box by $h(i, j)$. The \emph{hook multiset}
$H({\lambda})$ of the partition $\lambda$ is defined to be the
multiset of hook lengths of $\lambda.$  The \emph{conjugation} of
the partition $\lambda$ is defined by $\lambda'=(\lambda'_1,
\lambda'_2, \ldots, \lambda'_r)$ where $r=\lambda_1$ and
$\lambda'_i= \#\{ 1\leq j\leq m: \lambda_j\geq i \}$ for $1\leq i
\leq r.$ It is easy to see that, we can obtain the Young diagram of
the conjugation of a partition by reflecting its Young diagram along
its main diagonal. For example, Figure $1$ shows the Young diagrams
and hook lengths of partitions $(4,2,2)$ and $(3,3,1,1)$. It is easy
to see that these two partitions are conjugate to each other.

\begin{figure}[htbp]
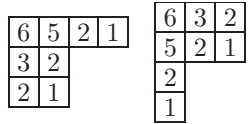

\begin{center}
\Yvcentermath1
\begin{tabular}{@{}cc@{}c@{}}
$\young(6521,32,21)$ & $\young(632,521,2,1)$
\end{tabular}

\end{center}
\caption{The Young diagrams and hook lengths of partitions $(4,2,2)$
and $(3,3,1,1)$.}
\end{figure}

The \emph{$\beta$-set} $\beta(\lambda)$ of the partition $\lambda$
is the set of hook lengths of boxes in the first column of the
corresponding Young diagram, which means that
$$\beta(\lambda)=\{\lambda_i+m-i : 1 \leq i \leq m\}.$$ It is easy
to see that a partition is uniquely determined by its $\beta$-set.
If $A$ is a finite set of some positive integers, we define
$\lambda_A$ to be the partition whose $\beta$-set is $A.$

 As we know, two conjugate partitions have the same hook multiset.
But in general, a partition is not uniquely determined by its hook
multiset up to conjugation. For instance, Herman and Chung
\cite{Herman} showed that for any $n\in \textbf{N}$, the following
two non-conjugate partitions
 have the same hook multiset:
$$\lambda_n = (n + 6, n + 3, n + 3, 2) \ \text{and}\ {\mu}_n = (n + 5, n + 5, n + 2, 1,
1).$$

Herman and Chung  proposed the following question:

\begin{ques} (\cite{Herman}.)
 What are the
necessary and sufficient conditions for partitions to be determined
by their hook multisets up to conjugation?
\end{ques}

\vspace{1ex}

Equivalently, given two different finite sets $A$ and $B$ of some
positive integers  with $H({\lambda}_A)=H({\lambda}_B),$ what are
the necessary and sufficient conditions for partitions ${\lambda}_A$
and ${\lambda}_B$ to be conjugate to each other?

We will give an answer to this question:

\begin{thm} \label{main}
Let $A$ and $B$ be  two different subsets of $\{ 1,2, \ldots, n \}$
with $n\in A\bigcap B$ and $H({\lambda}_A)=H({\lambda}_B).$ Then
partitions ${\lambda}_A$ and ${\lambda}_B$ are conjugate to each
other if and only if $A\setminus B$ and $B\setminus A$ are
$n-$symmetric sets, where a set $S$ is called  $n-$symmetric if
$S=\{n-x: x\in S \}.$
\end{thm}

\noindent\textbf{Example.} Let $A=\{2,3,6  \}$ and   $B=\{1,2,5,6
\}.$ In this case ${\lambda}_A=(4,2,2)$ and ${\lambda}_B=(3,3,1,1).$
By Figure $1$ we know $H({\lambda}_A)=H({\lambda}_B).$  Since
$A\setminus B=\{ 3\}$ and $B\setminus A=\{1,5 \}$ are $6-$symmetric
sets, by Theorem \ref{main} we know ${\lambda}_A=(4,2,2)$ and
${\lambda}_B=(3,3,1,1)$ must be conjugate to each other, which is
indeed true.

\section{Main results}

Suppose that $n$ is a positive integer and $A$ is a finite  set of
some positive integers. Let  $I_n=\{0,1, \ldots , n\}$ and
$n-A=\{n-x: x\in A \}$.  The following results are well-known and
easy to prove:

\begin{lem} \label{lem1}
(\cite{berge}.) Suppose that $A$ is a finite set of  some positive
integers whose largest element is $n$. Let $A'$ be the complement of
$A$ in $I_n$. Then the $\beta$-set of the conjugation of $\lambda_A$
is $n-A'.$ The hook multiset of $\lambda_A$ is
$$H(\lambda_A) = \{a - a' : a\in A, \ a' \in A' , a
> a'\}.$$
 \end{lem}

Now we can give the proof of Theorem \ref{main}.

\noindent\textbf{Proof of Theorem \ref{main}.} $\Rightarrow$:
Suppose that partitions ${\lambda}_A$ and ${\lambda}_B$ are
conjugate to each other. Let $A'$ be the complement of $A$ in $I_n$.
By Lemma \ref{lem1}, we know $B=n-A'$, thus $(n-A)\bigcup B= I_n$
and $(n-A)\bigcap B= \emptyset$. Let $a\in A\setminus B$. By $a\in
A$ we know $n-a\ \in n-A$, which means that $n-a \notin B$. By $a
\notin B$ we know $a\in n-A$, which means that $n-a\in A$. Then we
have $n-a\in A\setminus B$. Thus we know $A\setminus B$ is an
$n-$symmetric set. Similarly, $B\setminus A$ is also an
$n-$symmetric set.

$\Leftarrow$: Suppose that $A\setminus B$ and $B\setminus A$ are
$n-$symmetric sets. Let $A'$ and $B'$ be the complements of $A$ and
$B$ in $I_n$. We assume that ${\lambda}_A$ and ${\lambda}_B$ are not
conjugate to each other. Next we will deduce a contradiction to show
that this assumption couldn't be true.

\textbf{Step 1}: In this step we  will prove that $(n-A')\setminus
B$ and $B\setminus (n-A')$ are $n-$symmetric sets:

First let $a\in (n-A')\setminus B$.

Then we know $a\in (n-A')=n-I_n\setminus A=I_n\setminus (n-A)$. Thus
we know $$a\ \notin n-A=(n-(A\bigcap B))\bigcup (n-(A\setminus
B))=(n-(A\bigcap B))\bigcup (A\setminus B)$$ since $A\setminus B$ is
an $n-$symmetric set, which means that $a \notin A\setminus B$. But
we already know $a\notin B$. This implies that $a\notin A$, which
means that $a\in A' $ and thus $n-a\in n-A'$.

On the other hand, by $a\notin n-A$ we know $a\notin n-A\bigcap B$.
By $a\notin B$ we know $a\notin B\setminus A=n-B\setminus A$ since
$B\setminus A$ is an $n-$symmetric set. Put these two results
together, we have $$a\notin (n-A\bigcap B) \bigcup (n-B\setminus A)=
n-((A\bigcap B) \bigcup (n-B\setminus A))=n-B,$$ which means that
$n-a\notin B$.

Now we  have $n-a\in (n-A')\setminus B$ and thus $(n-A')\setminus B$
is an $n-$symmetric set.

Next let $b\in B\setminus (n-A')$.

Then we know $b\in B \bigcap (n-A).$ Thus $n-b\in A$, which means
that $n-b\notin B\setminus A=n-B\setminus A$ since $B\setminus A$ is
an $n-$symmetric set. Then we have  $b\notin B\setminus A$. But we
already know $b\in B$. Thus $b\in A$, which means that $n-b\notin
n-A'.$

On the other hand, $b\in B$ implies that $b\notin A\setminus
B=n-A\setminus B$  since $A\setminus B$ is an $n-$symmetric set.
Then we have $n-b\notin A\setminus B$. But we already know $n-b\in
A$, thus we know $n-b\in B$.

Now we have $n-b\in B\setminus (n-A')$ and thus $B\setminus (n-A')$
is an $n-$symmetric set.

\textbf{Step 2}: Suppose that $k\in \{1,2, \ldots, n\}$. Let $M$ and
$N$ be  multisets. We define $M_k$ to be the multiplicity of $k$ in
$M$. We also define $M+N$ to be the multiset $\{x+y:x\in M, \ y\in N
\}$ and $M-N$ to be the multiset $\{x-y:x\in M, \ y\in N \}$.
 By Lemma \ref{lem1},
the multiplicity of $k$ in $H({\lambda}_A)$ and $H({\lambda}_B)$ are
${H({\lambda}_A)}_k=(A-A')_k$ and ${H({\lambda}_B)}_k=(B-B')_k$. We
know ${H({\lambda}_A)}_k={H({\lambda}_B)}_k$ since
$H({\lambda}_A)=H({\lambda}_B).$ Let $C=A\setminus B,\  D=B\setminus
A$ and $E=A\bigcap B.$ We know $C=n-C,\ D=n-D$ since $A\setminus B$
and $B\setminus A$ are $n-$symmetric sets.

 We mention that, all the following set operations in  Step 2 are on
multisets. Then we have

\begin{eqnarray*}
{H({\lambda}_A)}_k-{H({\lambda}_B)}_k  & = &( A-A')_k-(B-B')_k \\
&= & (C\bigcup E-A')_k-(D\bigcup E-B')_k  \\  & = &(C-A')_k+(E-A')_k
 -  ( D -B')_k- (E -B')_k.
\end{eqnarray*}

Since $C=n-C,\ D=n-D$, we have
\begin{eqnarray*}
(E-A')_k - (E -B')_k & =  & (E-I_n\setminus (C\bigcup E))_k- (E
-I_n\setminus (D\bigcup E))_k \\ &= &-(E-C)_k+(E-D)_k \\ & = &
-(E+n-C)_{n+k}+(E+n-D)_{n+k}\\ & = & -(E+C)_{n+k}+(E+D)_{n+k}.
\end{eqnarray*}

Then we know

\begin{eqnarray*}
& &{H({\lambda}_A)}_k-{H({\lambda}_B)}_k \\ & =
&(C-A')_k-(E+C)_{n+k}
- ( D -B')_k+(E+D)_{n+k} \\ & = &(C-A')_k-(E+C)_{n+k}-(D+C)_{n+k} \\
& + &(D+n-C)_{n+k}- ( D -I_n\setminus B)_k+(E+D)_{n+k}
\\ & = &(C+n-A')_{n+k}-(B+C)_{n+k} \\ & +
&(D+n-C)_{n+k}- ( D+n -I_n\setminus B)_{n+k}+(E+D)_{n+k}
\\ & = &(C+n-A')_{n+k}-(B+C)_{n+k} \\ & +
&(D+n-C)_{n+k}- ( D+ I_n\setminus (n-B))_{n+k}+(E+D)_{n+k}
\\ & = &(C+n-A')_{n+k}-(B+C)_{n+k}   +
(D+n-C)_{n+k} \\ &-& ( D+ I_n\setminus (n-D))_{n+k}+( D+
(n-E))_{n+k}+(E+D)_{n+k}
\\ & = &(C+n-A')_{n+k}-(B+C)_{n+k} \\ & +
&(D+n-C\bigcup E)_{n+k}- ( D+ I_n\setminus D)_{n+k}+(E+D)_{n+k}
\\ & = &(C+n-A')_{n+k}-(B+C)_{n+k} \\ & +
&(D+n-A)_{n+k}- ( D+ I_n)_{n+k}+(D+ B)_{n+k}
\\ & = &(C+n-A')_{n+k}-(B+C)_{n+k} \\ & -
&( D+ I_n\setminus (n-A))_{n+k}+(D+ B)_{n+k}
\\  & =
&(C+n-A')_{n+k}-(B+C)_{n+k} \\&-& ( D+n -A')_{n+k}+(D+B)_{n+k}.
\end{eqnarray*}

By Lemma \ref{lem1} we know $C\bigcup D$ and $((n-A')\setminus
B)\bigcup (B\setminus (n-A') )$ are not empty sets since we already
assume that ${\lambda}_A$ and ${\lambda}_B$ are not conjugate to
each other.

 Let $x$ be the greatest element in $C\bigcup D$ and $y$ be the
greatest element in $((n-A')\setminus B)\bigcup (B\setminus (n-A')
)$.  We have $x\geq \frac{n}{2}$ since $C$ and $D$ are $n-$symmetric
sets. By Step 1, $(n-A')\setminus B$ and $B\setminus (n-A') $ are
also $n-$symmetric sets, then we have $y\geq \frac{n}{2}$. It is
easy to see that, $x= \frac{n}{2}$ and $y= \frac{n}{2}$ couldn't be
true simultaneously. Let $z=x+y-n$. Thus we have $1\leq z \leq n$.

Then one of the following cases must be true:

(1) $x\in C, y\in (n-A')\setminus B$;

(2) $x\in C, y\in B\setminus (n-A')$;

(3) $x\in D, y\in (n-A')\setminus B$;

(4) $x\in D, y\in B\setminus (n-A')$.

For case $(1)$, if $x\in C, y\in (n-A')\setminus B$, then it is easy
to see that $$(C+n-A')_{n+z}\geq 1$$ and $$(B+C)_{n+z} = ( D+n
-A')_{n+z}=(D+B)_{n+z}=0$$ since $x+y=n+z$, $x$ is the greatest
element in $C\bigcup D$ and $y$ is the greatest element in
$((n-A')\setminus B)\bigcup (B\setminus (n-A') )$. Then we have
\begin{eqnarray*}  && {H({\lambda}_A)}_z-{H({\lambda}_B)}_z\\ &=& (C+n-A')_{n+z} -(B+C)_{n+z}
-( D+n -A')_{n+z}+(D+B)_{n+z} \geq 1, \end{eqnarray*} which means
that $H({\lambda}_A)\neq H({\lambda}_B),$ a contradiction!

For case $(2)$, $(3)$ or $(4)$, similarly we can also deduce that
$H({\lambda}_A)\neq H({\lambda}_B),$ which is a contradiction.

By this contradiction we know ${\lambda}_A$ and ${\lambda}_B$ must
be conjugate to each other. We finish the proof. \hfill $\square$

% ------------------------------------------------------------------------

\section{Acknowledgements}
The author  thanks Prof. P. O. Dehaye for the encouragement and
helpful discussions. The author is supported  by Forschungskredit of
the University of Zurich, grant no. [FK-14-093].

\end{document}